\documentclass[preprint,12pt]{article}
\usepackage[latin9]{inputenc}
\usepackage{float}
\usepackage{amsmath}
\usepackage{amssymb}
\usepackage{graphicx}
\usepackage{esint}

\usepackage{amsmath,amstext,amsfonts,amssymb,amsbsy,array,graphicx}
\usepackage{stmaryrd}
\usepackage{float}
\usepackage{subfig}

\usepackage[english]{babel}

\def\IR{\mathbb R}
\newcommand{\bfI}{\boldsymbol I}

\newcommand{\bft}{\boldsymbol t}
\newcommand{\bfg}{\boldsymbol g}
\newcommand{\bfPs}{{\boldsymbol P}_\Gamma}

\newcommand{\bfx}{{\boldsymbol x}}

\newcommand{\mcK}{\mathcal{K}}
\newcommand{\mcF}{\mathcal{F}}
\newcommand{\bfn}{{\boldsymbol n}}

\newcommand{\bfu}{{\boldsymbol u}}
\newcommand{\bfv}{{\boldsymbol v}}

\newcommand{\bfw}{{\boldsymbol w}}

\usepackage{subfig}
\makeatother

\begin{document}
\title{Minimal surface computation using a finite element method on an embedded surface}
\author{Mirza Cenanovic$^a$, Peter Hansbo$^a$, Mats G. Larson$^b$\\[4mm] \it\small$^a$Department of Mechanical Engineering, \\\it\small J\"onk\"oping University\\\it\small
SE-55111 J\"onk\"oping, Sweden \\\it\small $^b$Department of Mathematics and Mathematical Statistics\\\it\small Ume{\aa} University\\\it\small 
SE-901 87 Ume{\aa}, Sweden}
\date{}
\maketitle
\begin{abstract}
We suggest a finite element method for computing minimal surfaces based
on computing a discrete Laplace--Beltrami operator operating on the coordinates of the surface.
The surface is a discrete representation of the zero level set of a distance function using linear tetrahedral finite elements, 
and the finite element discretization is done on the piecewise planar isosurface using the shape functions
from the background three dimensional mesh used to represent the distance function. A recently suggested stabilization scheme
is a crucial component in the method.
\end{abstract}
%

\section{Introduction}

An important application of partial differential equations on surfaces
is equations that are used to determine the shape of surfaces in order
to satisfy certain design criteria, called \emph{form finding}, cf.
\cite{BoKoPaAlLe10}. A classical example is the minimal surface problem,
the simplest kind of form finding, where a surface with minimal curvature
is sought, given the position of its boundary.

In this paper, we consider minimizing curvature by means of discretizing
the Laplace-Beltrami operator on an embedded surface, following Olshanskii,
Reusken, and Grande \cite{OlReGr09}, where a Galerkin method is
constructed by using the restrictions of continuous piecewise linears
defined on embedding tetrahedra to the embedded surface. The Laplace--Beltrami
operator applied to the Cartesian coordinates of the surface gives
the curvature vector $(\kappa_1+\kappa_2)\bfn$, where $\kappa_1$ and $\kappa_2$ are the the principal curvatures and $\bfn$ is the surface normal,
cf. \cite{BoKoPaAlLe10}. Our algorithm is motivated by previous work on meshed
surfaces by Dziuk \cite{Dz91}, where viscous relaxation was used
to move a triangulated surface on which the Laplace--Beltrami operator
was discretized. A method on embedded surfaces related to ours is
given by Chopp \cite{Ch92}, where the curvature is computed using
instead the values of a level set implicitly defining the embedded
surface. Classically, this problem has also been considered in a referential
domain, cf. Johnson and Thom\'ee \cite{JoTh75}, which severely limits
the scope of the method.

\section{Model problem and finite element method}

\subsection{The continuous problem}

Let $\Gamma(t)$ denote a time--dependent smooth two-dimensional surface
embedded in $\IR^{3}$ with signed distance function $\varrho$. For
notational convenience, we shall frequently omit the dependence on
$t$ in $\Gamma(t)$.

We consider the following problem: given a final time $T$, find $\bfx_{\Gamma}:\Gamma(t)\rightarrow\IR^{3}$
such that 
\begin{align}
\dot{\bfx}_{\Gamma}-\Delta_{\Gamma(t)}\bfx_{\Gamma}=0\quad\text{in \ensuremath{\Gamma(t),\; t\in(0,T),}}\label{eq:LB}
\end{align}
${\bfx}_{\Gamma}=\bfg$ on $\partial\Gamma(t)$ (we shall assume $\bfg$
is constant). Here 
\[
\dot{\bfx}_{\Gamma}:=\frac{\partial\bfx_{\Gamma}}{\partial t},
\]
$\Delta_{\Gamma}$ is the Laplace-Beltrami operator defined by 
\begin{equation}
\Delta_{\Gamma}=\nabla_{\Gamma}\cdot\nabla_{\Gamma}
\end{equation}
where $\nabla_{\Gamma}$ is the surface gradient 
\begin{equation}
\nabla_{\Gamma}=\bfPs\nabla
\end{equation}
with $\bfPs=\bfPs(\bfx)$ the projection of $\IR^{3}$ onto the tangent
plane of $\Gamma$ at $\bfx\in\Gamma$, defined by 
\begin{equation}
\bfPs=\bfI-\bfn\otimes\bfn
\end{equation}
where $\bfn=\nabla\varrho$ denotes the exterior normal to $\Gamma$
at $\bfx$ , $\bfI$ is the identity matrix.

As is well known, cf. \cite{Dz91}, 
\[
-\Delta_{\Gamma(t)}\bfx_{\Gamma}=2H\bfn
\]
where $H=(\kappa_1+\kappa_2)/2$ is the mean curvature of $\Gamma(t)$. Thus, we have the
curvature driven normal flow 
\[
\dot{\bfx}_{\Gamma}=-2H\bfn.
\]
If ${\bfx}_{\Gamma}$ is found by following the zero isosurface of a
level set function $\phi(\bfx,t)$, then the material derivative of the level set function at ${\bfx}_{\Gamma}$ is given by
\begin{equation}
\frac{d\phi}{dt}=\frac{\partial\phi}{\partial t}+\dot{\bfx}_{\Gamma}\cdot\nabla\phi=0,\label{dist}
\end{equation}
and, since $\nabla\phi=\bfn$, if we assume that $\vert\nabla\phi\vert=1$
is enforced (i.e., the level set function is a distance function), we may update the level set function
at the surface directly via 
\[
\frac{\partial\phi}{\partial t}=2H.
\]
This idea, together with a relation between $H$ and spatial derivatives
of $\phi$, was used by Chopp \cite{Ch92}. Here we shall instead
use (\ref{eq:LB}) directly.

Defining 
\[
V=\{\bfv\in[H^{1}(\Gamma)]^{3}:\;\bfv={\bold0}\;\text{on}\;\partial\Gamma\},
\]
the weak statement corresponding to (\ref{eq:LB}) takes the form:
given the coordinate map or embedding of $\Gamma=\Gamma(t)$ into
$\IR^{3}$, denoted by $\bfx_{\Gamma}:\Gamma\ni\bfx\mapsto\bfx\in\IR^{3}$,
find the velocity $\dot{\bfx}_{\Gamma} =: \bfu_{\Gamma}\in V$ such that 
\begin{equation}
({\bfu}_{\Gamma},\bfv)_{\Gamma(t)}+a(\bfx_{\Gamma},\bfv)=0\quad\forall\bfv\in V\label{cont}
\end{equation}
where 
\begin{equation}
a(\bfu,\bfv)=(\nabla_{\Gamma}\bfu,\nabla_{\Gamma}\bfv)_{\Gamma(t)},
\end{equation}
and $(\bfv,\bfw)_{\Gamma}=\int_{\Gamma}\bfv\cdot\bfw\, d\Gamma$ is the $L^{2}$
inner product on $\Gamma$.

\subsection{Discretization in time and space}

\label{Sec:FEM}

In Dziuk \cite{Dz91}, a semi--discrete version of (\ref{cont})
of backward Euler type on meshed surfaces, yielding a discrete velocity $\bfu_{\Gamma}^h$, so that given the nodal
coordinates ${\bold x}_{N}^{n}$, ${\bold x}_{N}^{n}\approx{\bold x}_{N}(t_{n})$,
${\bold x}_{N}^{n+1}$ was computed by 
\begin{equation}
{\bold x}_{N}^{n+1}={\bold x}_{N}^{n}+k_{n}\bfu_{\Gamma}({\bold x}_{N}^{n}),\label{update}
\end{equation}
where $k_{n}=t_{n+1}-t_{n}$, updating the mesh, updating $n$, and
continuing until the curvature is small enough. We now wish to instead
solve this problem on embedded surfaces in $\IR^{3}$ using the general
technique proposed by Olshanskii et al. \cite{OlReGr09} by means
of a distance function for the definition of $\Gamma$ and we cannot
thus directly solve for the location of $\Gamma$ but need to solve
the additional equation (\ref{dist}) for the level set function.

\begin{figure}[H]
\centering{}\includegraphics[width=0.5\paperwidth]{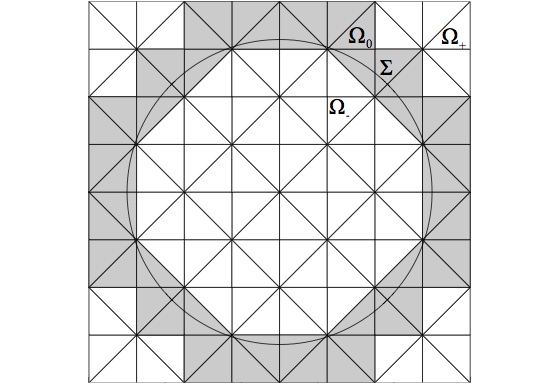}\caption{Intersecting domains\label{fig:Intersecting-domains}}
\end{figure}

To discretize in space, let $\mcK$ be a quasi uniform partition into
shape regular tetrahedra of a domain $\Omega$ in $\IR^{3}$ completely
containing $\Gamma^{n}$ for all $n$. Furthermore, we let 
\begin{equation}
\mcK_{h}^{n}=\{K\in\mcK:K\cap\Gamma^{n}\neq\emptyset\}\label{eq:cutEle}
\end{equation}
be the set of tetrahedra that intersect $\Gamma^{n}$. See Figure
\ref{fig:Intersecting-domains}. 

\begin{figure}[H]
\begin{centering}
\includegraphics[width=0.5\paperwidth]{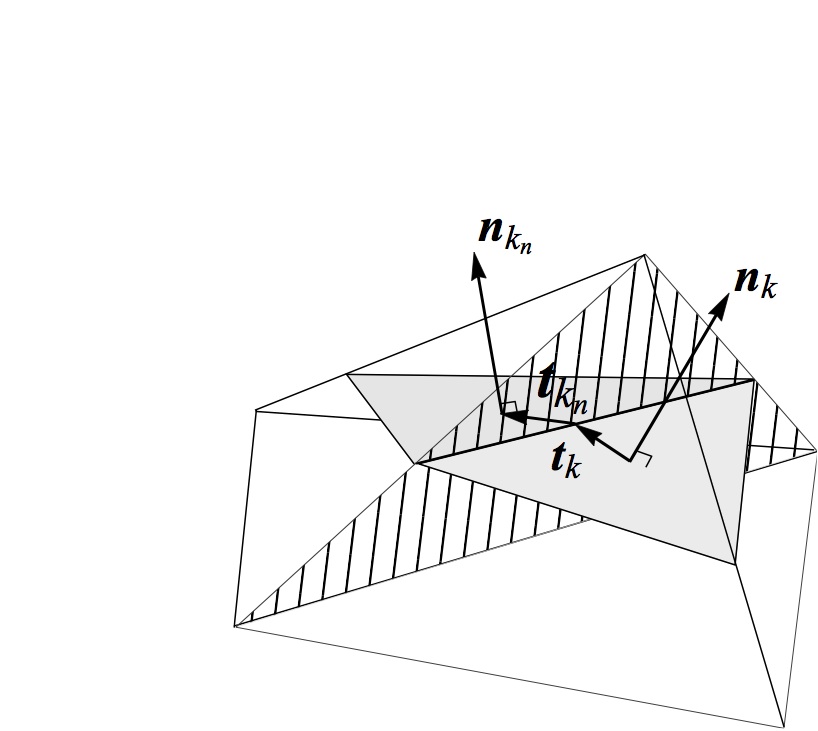}
\par\end{centering}

\centering{}\caption{Tangential jump\label{fig:Tangential-jump}}
\end{figure}

We let 
\[
V_{h}^{n}=\{\bfv:\;\text{$\bfv$ is a continuous piecewise linear polynomial defined on \ensuremath{\mcK_{h}^{n}}},\;\bfv={\bold0}\;\text{on}\;\partial\Gamma^{n}\}.
\]
In general, a finite element method for computing the mean curvature vector
by means of the discrete Laplacian may fail due to instability, see
\cite{HaLaZa13} where the following stabilization method was suggested:
given $\Gamma^{n}$ and the corresponding coordinate function $\bfx_{\Gamma}^{n}$, find
$\bfu_{\Gamma}^{h}\in V_{h}^{n}$ such that 
\begin{equation}\label{contstab}
({\bfu}_{\Gamma}^{h},\bfv)_{\Gamma^{n}}+j({\bfu}_{\Gamma}^{h},\bfv)=-a(\bfx_{\Gamma}^{n},\bfv)_{\Gamma^{n}}\quad\forall\bfv\in V_{h}^{n}
\end{equation}
where the bilinear form $j(\cdot,\cdot)$ is defined by 
\begin{equation}
j(\bfu,\bfv)=\sum_{F\in\mcF_{I}}([\bft_{F}\cdot\nabla\bfu],[\bft_{F}\cdot\nabla\bfv])_{F}.
\end{equation}
Here, $\mcF_{I}$ denotes the set of internal interfaces in $\mcK_{h}^{n}$,
$[\bft_{F}\cdot\nabla v]=(\bft_{F}\cdot\nabla v)^{+}-(\bft_{F}\cdot\nabla v)^{-}$
with $w(\bfx)^{\pm}=\lim_{s\rightarrow0^{+}}w(\bfx\mp s\bft_{F})$,
is the jump in the tangent gradient across the face $F$. Here, the
jump in the tangent gradient at an edge $E$ shared by the elements
$K_{1}$ and $K_{2}$ is defined by 
\begin{equation}
[\bft_{E}\cdot\nabla\bfu]=\bft_{E,K_{1}}\cdot\nabla\bfu_{1}+\bft_{E,K_{2}}\cdot\nabla\bfu_{2}
\end{equation}
where $\bfu_{i}=\bfu|_{K_{i}}$, $i=1,2,$ and $\bft_{E,K_{i}}$ denotes
the outwards unit vector orthogonal to $E$, tangent to $K_{i}$,
$i=1,2.$ See figure \eqref{fig:Tangential-jump}
This stabilization method was shown to yield first order convergence of the curvature in $L_2(\Gamma)$
in the finite element method proposed in \cite{HaLaZa13}.

After obtaining ${\bfu}_{\Gamma}^{h}$ we use \eqref{dist} to get
the level set evolution equation 
\begin{equation}
\frac{\partial\phi}{\partial t}=\bfu_{\Gamma}^{h}\cdot\bfn=0,\label{eq:surf}
\end{equation}
where we employ a time discretization scheme to solve for the level set function at time $t_{n+1}$,  $\phi^{n+1}$. 

Note that this method requires that the level set function is a distance function, $|\nabla\phi|$=1. 
This is done by reinitialization
of the level set function. In order to keep computational effort
at a minimum both the reinitialization and the propagation are being
done on a narrow band of tetrahedral elements so that only a small
set of elements that are cut by the surface and their neighbor elements
are updated.

\section{Numerical implementation\label{implement}}

Below follows a step by step finite element implementation.

\begin{enumerate}
\item Construct a linear tetrahedral mesh $\mathcal{K}$ is created on the domain
$\Omega$ in $\mathbb{R^{\mathrm{3}}}$ in which we embedd the the
implicit surface $\Gamma$. Let ${\bold x}_N^n$ denote the vector of 
node coordinates in $\mathcal{K}$. 
\item Set up the level set function $\phi(\bfx,t_n)$ such that $\phi(\bfx_\Gamma,t_n)=0$. 
\item Discretize the distance function $\phi^{h,n}(\bfx)\approx \phi(\bfx,t_n)$ by evaluating it in the nodes
of the tetrahedral mesh giving a nodal vector $\boldsymbol{\phi}_N^n=\phi({\bold x}_N,t_n)$.
\item Initialize ${\bold x}_{N}^{0}={\bold x}_{N}(0)$
\item Find elements that are cut by the surface, $\mcK_{h}^{n}$, using
\eqref{eq:cutEle}. 
\item Extract zero isosurface points giving $\Gamma_{h}^{n}$ by going over all elements
$\mcK_{h}^{n}$ interpolating the signed distance function linearly
using the tetrahedral basis functions.
\item Compute the velocity field
$\boldsymbol{u}_{\Gamma^n}^{h}$ by solving \eqref{contstab}. This is done by solving the matrix
equation
\begin{equation}
{\bold M}\,\mathbf{u}_{N} =-{\bold S}\,\mathbf{x}_{N}^{n}
\end{equation}
where ${\bold u}_N$ are the nodal velocities in the band of elements containing $\Gamma^n$, $\bold S$ is the matrix corresponding to the Laplace--Beltrami operator, and $\bold M$ a stabilized mass matrix computed on $\Gamma^n$.
\item Choose a time step $k_n$ and propagate ${\phi}^{h,n}$ to ${\phi}^{h,n+1}$, in the nodes of the band of elements containing $\Gamma^n$, using
\begin{equation}
{\phi}^{h,n+1}={\phi}^{h,n} - k_{n}\boldsymbol{n}\cdot\bfu^h_{\Gamma^n}
\end{equation}
where 
\begin{equation}
\boldsymbol{n}=\frac{\nabla\phi^{h,n}}{\vert\nabla\phi^{h,n}\vert}
\end{equation}
is the normal vector at $\bfx$ (in practice we also use an $L_{2}$ projection
to represent $\bfn$ at the nodes).
\item Reinitialize $\boldsymbol{\phi}^{h,n+1}$: we reinitialize the function
near the front by computing on a narrow band of elements using the time step
as a guide for how many elements to select. The reinitialization can
be done in several ways see, e.g., \cite{Se99}; we chose to reinitialize by updating
the value of ${\phi}^{h,n+1}$ in each node in
the narrow band with the signed distance to the closest node on the
discrete surface $\Gamma_{h}^{n+1}$.
\item Compute the discrete mean curvature 
\[
H^h := -\frac{\bfn\cdot\boldsymbol{u}_{\Gamma^n}^{h}}{2},
\]
which should converge to zero everywhere.
\item If the $L_2$ norm $\| H^h\|_{\Gamma_h}\leq\epsilon$, where $\epsilon$ is a small number,
then stop; otherwise update $n \mapsto n+1$ 
and repeat from step 5. 
\end{enumerate}

\section{Numerical examples}

The following Figures show some examples of minimal surfaces computed
using the scheme of Section \ref{implement}. We give some examples on different surface-evolutions. In all convergence plots the curvature has been computed as the Euclidean norm of $H^h({\bold x}_N)$ (mean curvature evaluated in the nodes).

\begin{figure}
\begin{centering}
\subfloat[Surface at timestep 1]{\includegraphics[width=0.25\paperwidth]{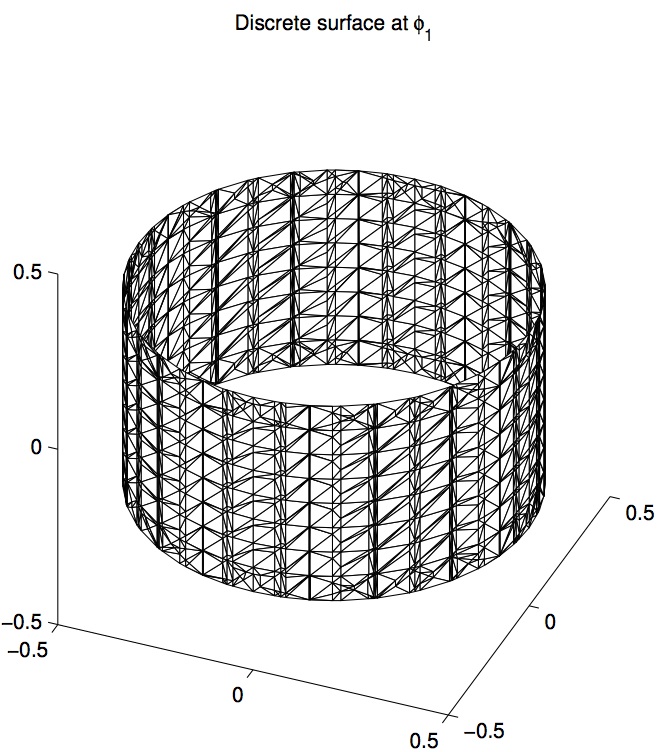}

}\subfloat[Surface at timestep 2]{\includegraphics[height=0.25\paperwidth]{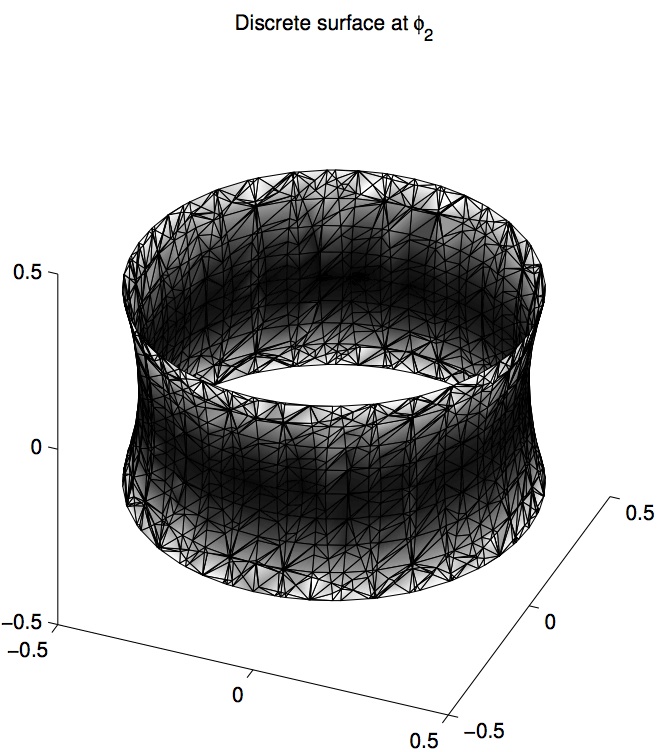}}
\par\end{centering}

\centering{}\subfloat[Surface at converged configuration]{\includegraphics[height=0.25\paperwidth]{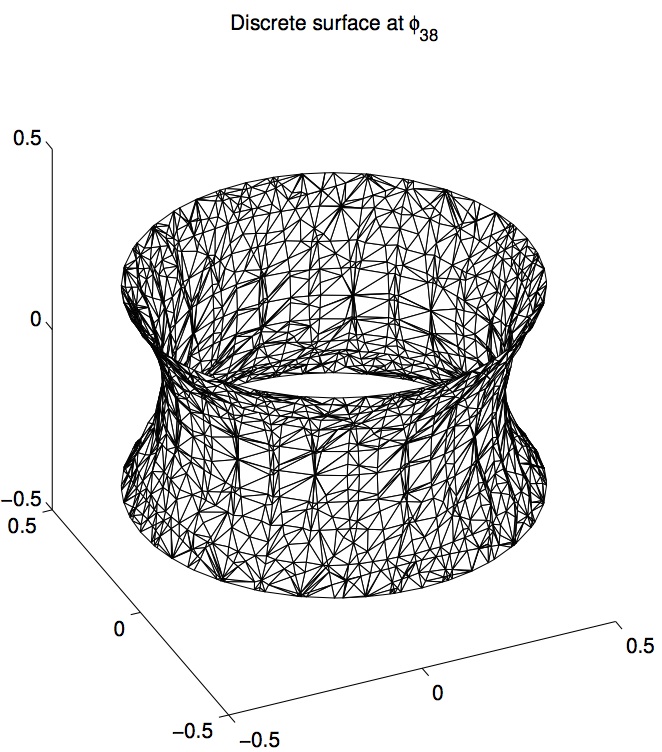}

}\subfloat[Top view of final configuration]{\includegraphics[height=0.25\paperwidth]{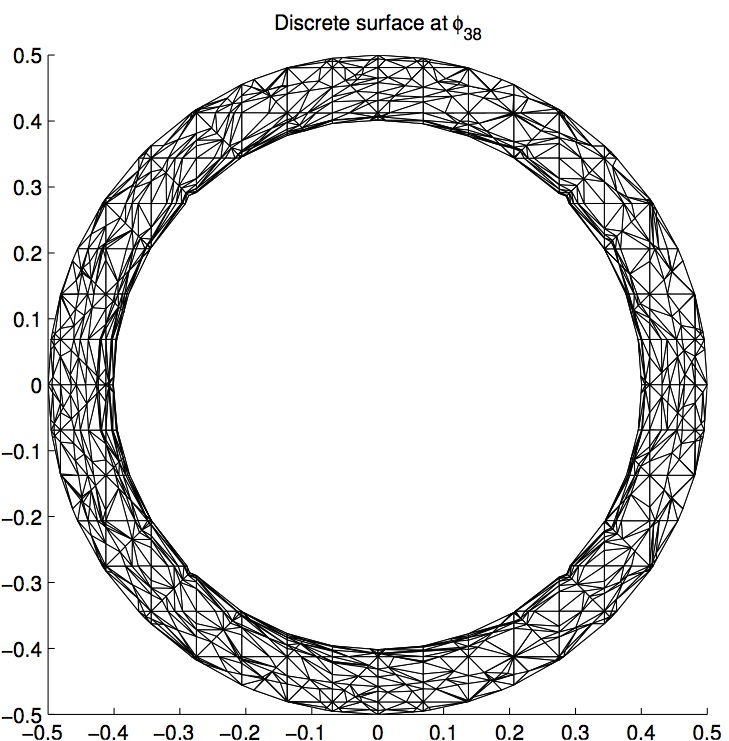}}\caption{Catenoid\label{fig:Catenoid}}
\end{figure}

\begin{figure}
\begin{centering}
\subfloat[Surface at timestep 1]{\includegraphics[width=0.25\paperwidth]{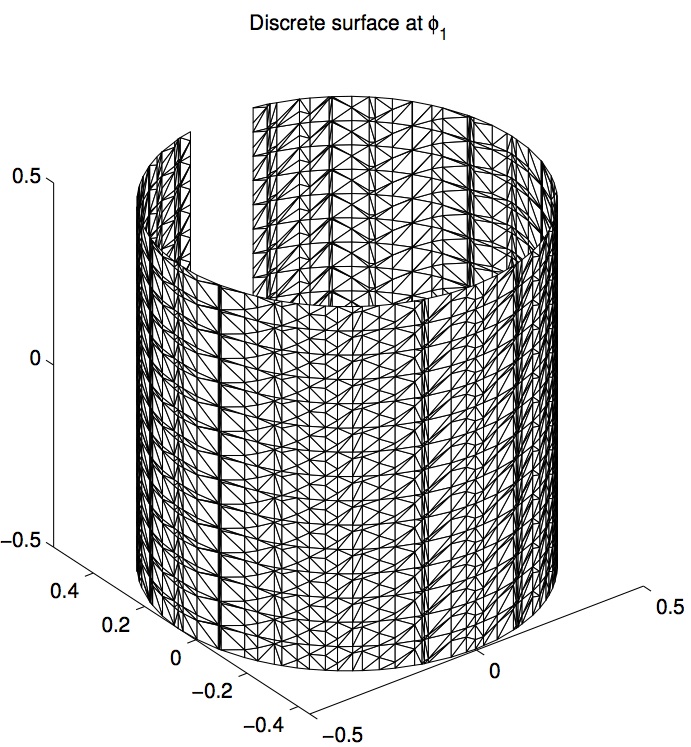}}\subfloat[Surface at timestep 2]{\includegraphics[height=0.25\paperwidth]{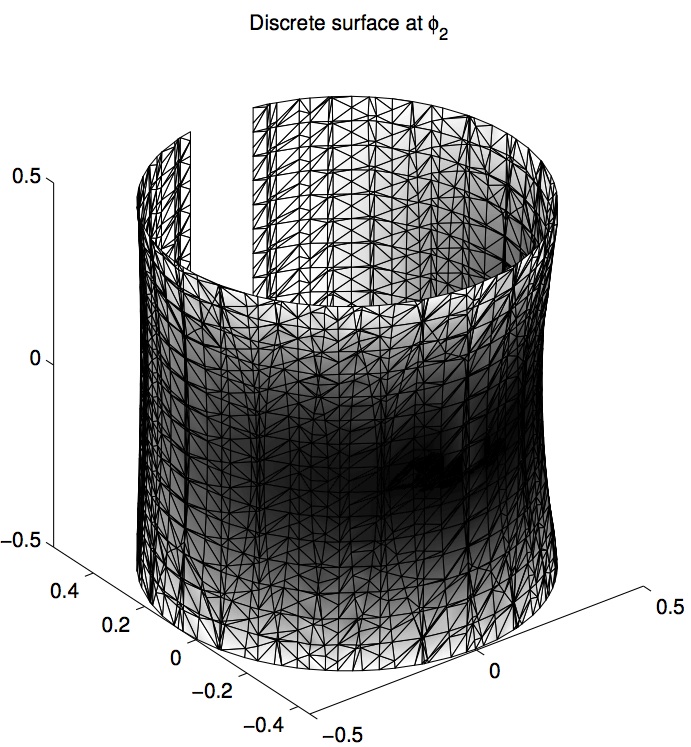}}
\par\end{centering}

\centering{}\subfloat[Surface at converged configuration]{\includegraphics[height=0.25\paperwidth]{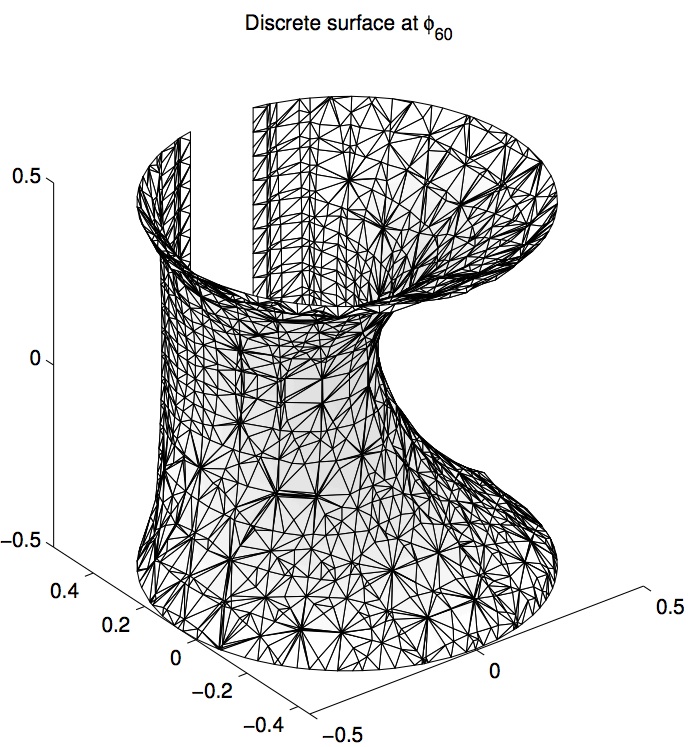}}\subfloat[Side view]{\includegraphics[height=0.25\paperwidth]{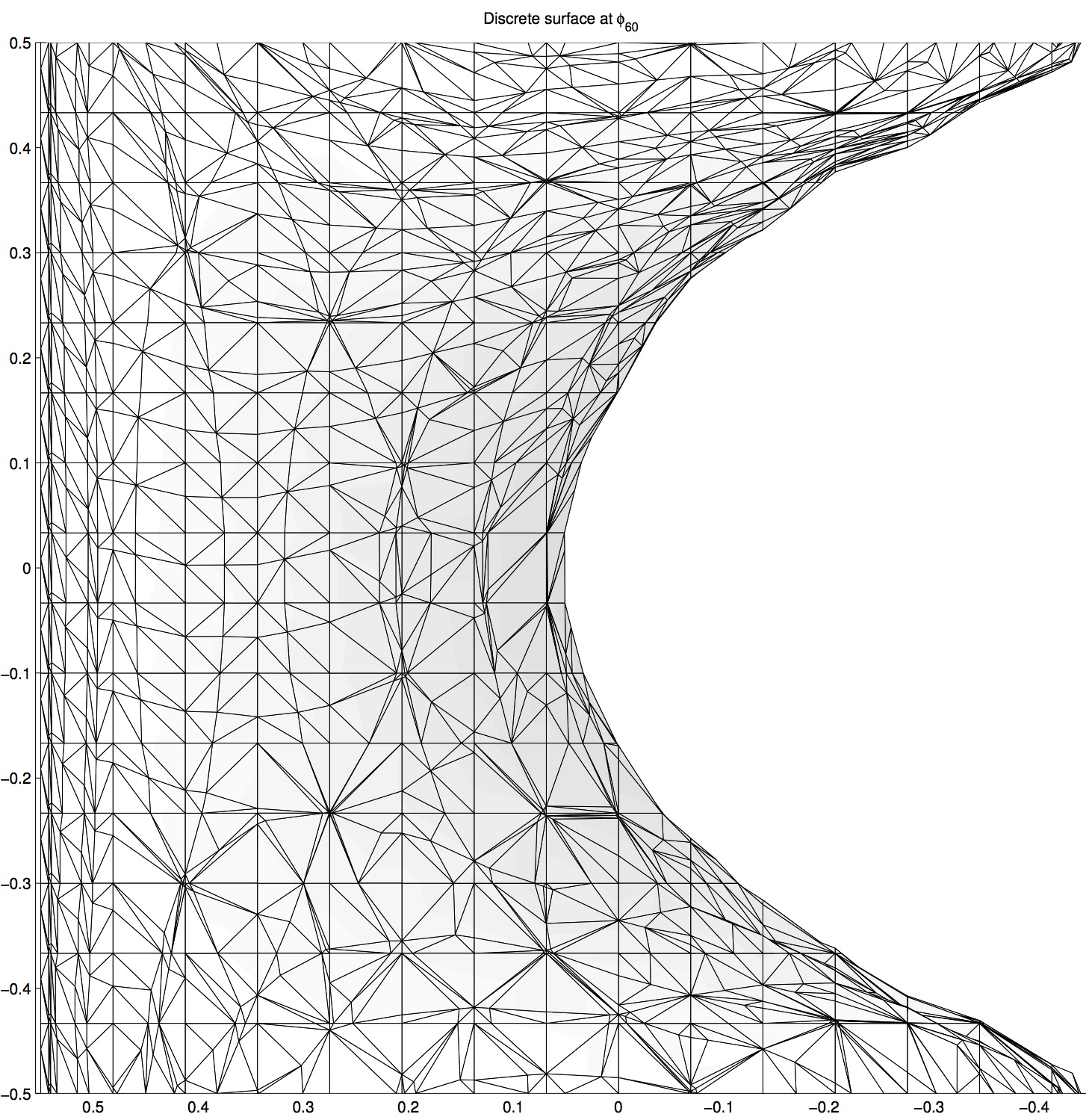}}\caption{Cut catenoid\label{fig:CutCatenoid}}
\end{figure}

Fig \ref{fig:Catenoid} shows a cylinder evolving to a catenoid. The
initial radius is 0.5, with axis centered at $x=0$, $y=0$, the dimensions are taken from \cite{Ch92}; the height is
0.554518 and the mesh domain is of the size $1.2\times 1.2\times 0.554518$
and consists of 14761 tetrahedral elements. Figure \ref{fig:Catenoid} shows
the evolutions from a cylinder to a catenoid with a inner diameter
of approximately $0.4$ which agrees with the theoretical result. The face colors represent a velocity field
where dark is high velocity. Figure \ref{fig:CutCatenoid} shows the
convergence of mean curvature for the catenoid.

\begin{figure}
\begin{centering}
\includegraphics[width=0.5\columnwidth]{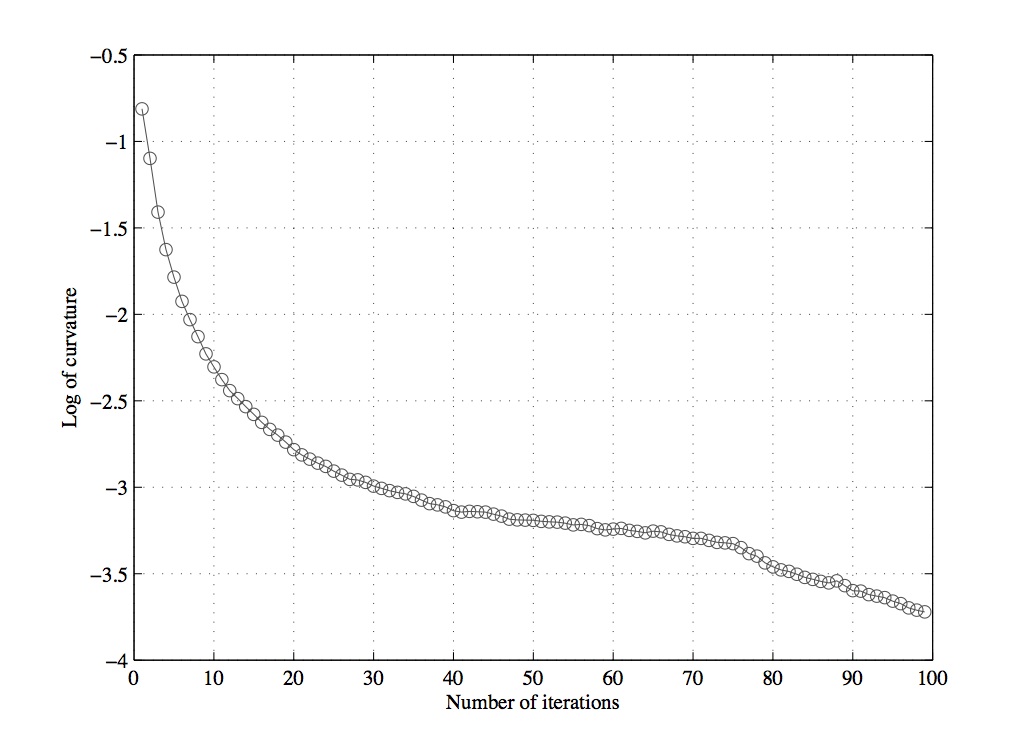}\caption{Convergence of the cut catenoid curvature\label{fig:CutCatenoidConverg}}

\par\end{centering}

\end{figure}

In Fig \ref{fig:CutCatenoid} a cut cylinder is shown. The initial
radius is 0.5, with axis positioned at $x=0$, $y=0.06$. The height
is 0.554518 and the mesh domain is of the size $1.2\times1.2\times0.554518$,
there are 23040 tetrahedral elelements. Figure \ref{fig:CutCatenoidConverg}
shows the evolution from a cut cylinder to a cut catenoid. Figure
\ref{fig:CutCatenoidConverg} shows the convergence of mean curvature for the cut catenoid.

\begin{figure}
\begin{centering}
\subfloat[Surface at timestep 1]{\includegraphics[width=0.25\paperwidth]{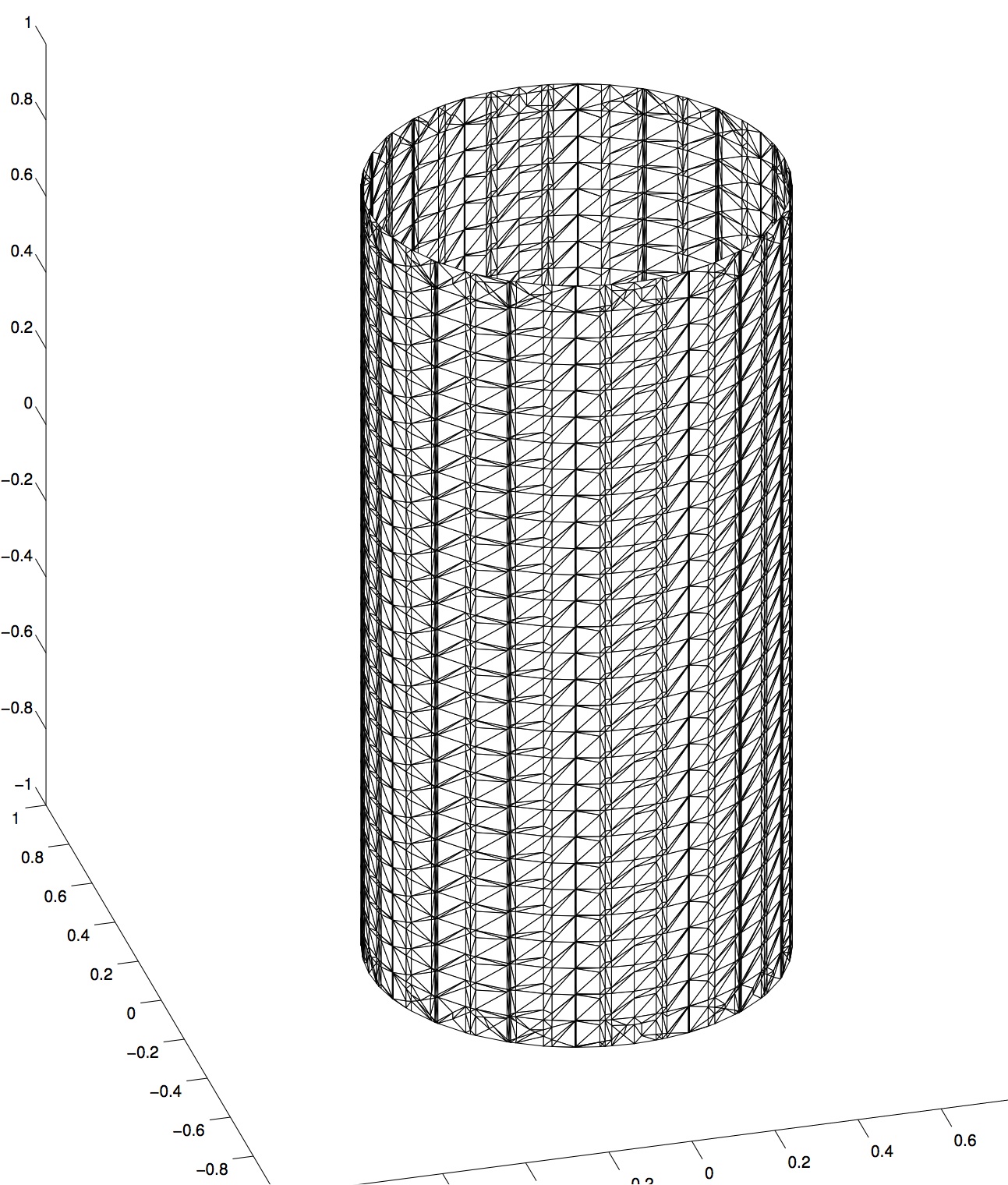}}
\subfloat[Surface at timestep 4]{\includegraphics[height=0.25\paperwidth]{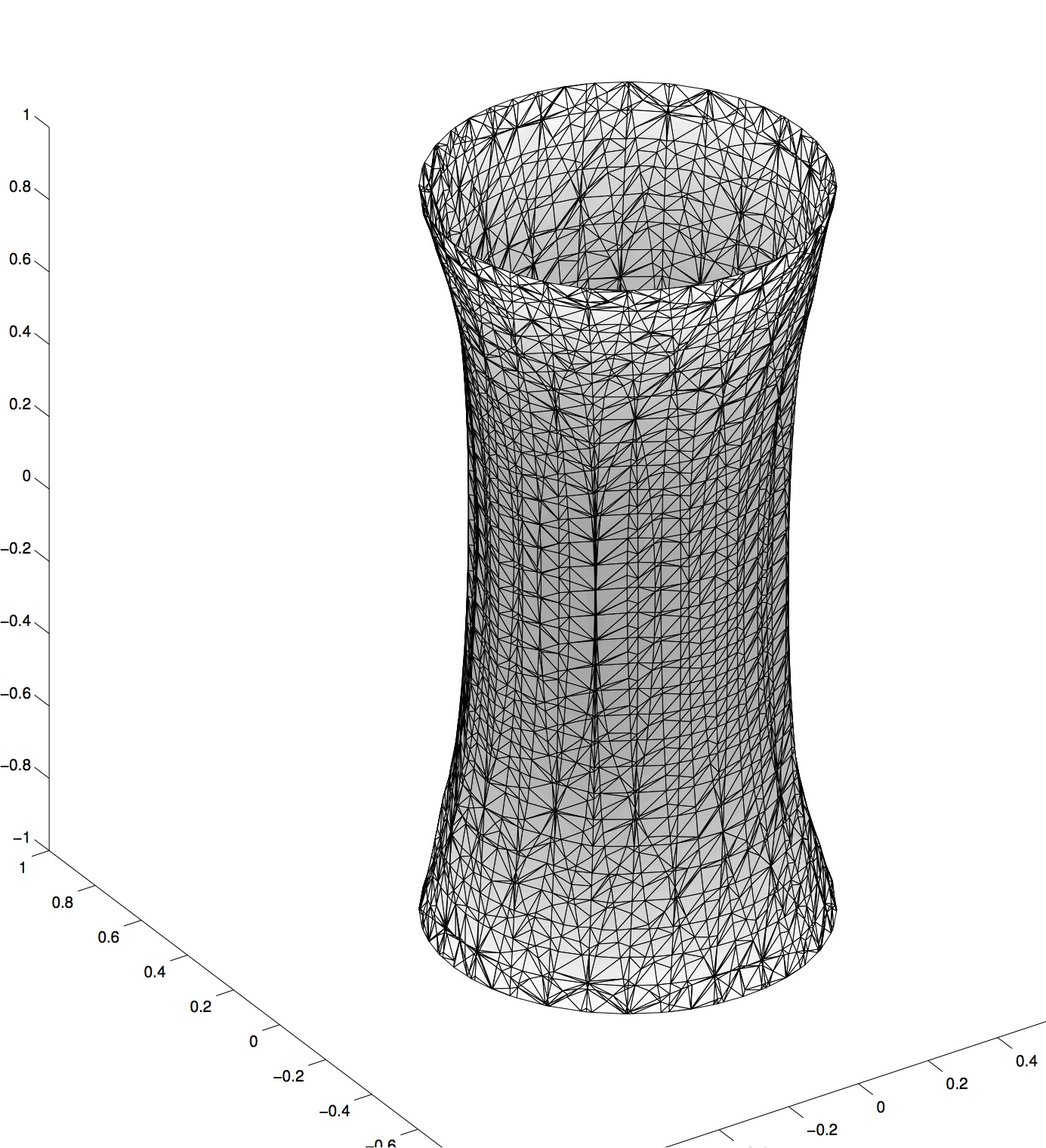}}
\par\end{centering}

\begin{centering}
\subfloat[Surface right before collapsing]{\includegraphics[height=0.25\paperwidth]{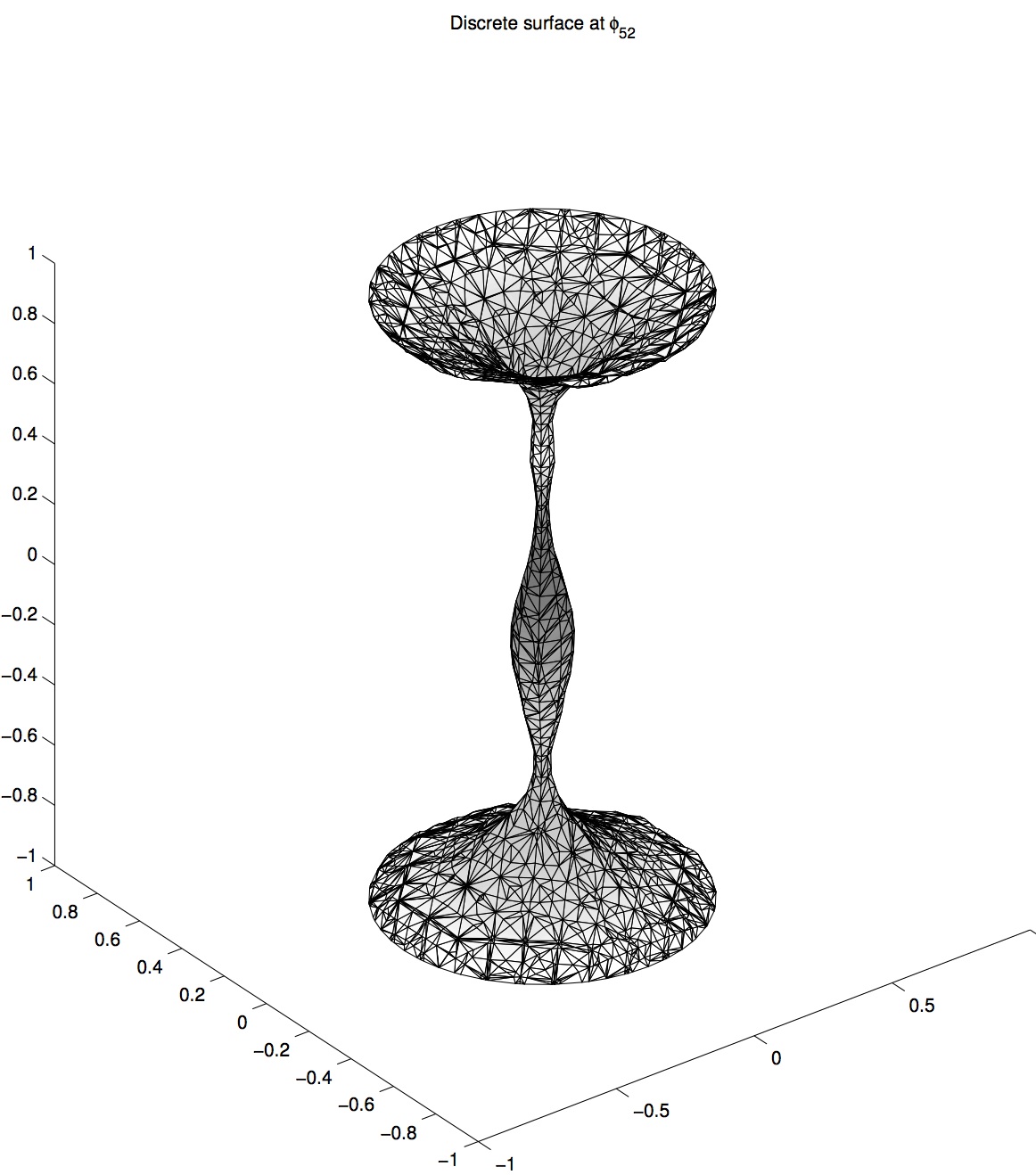}}
\subfloat[Collapsed surface]{\includegraphics[height=0.25\paperwidth]{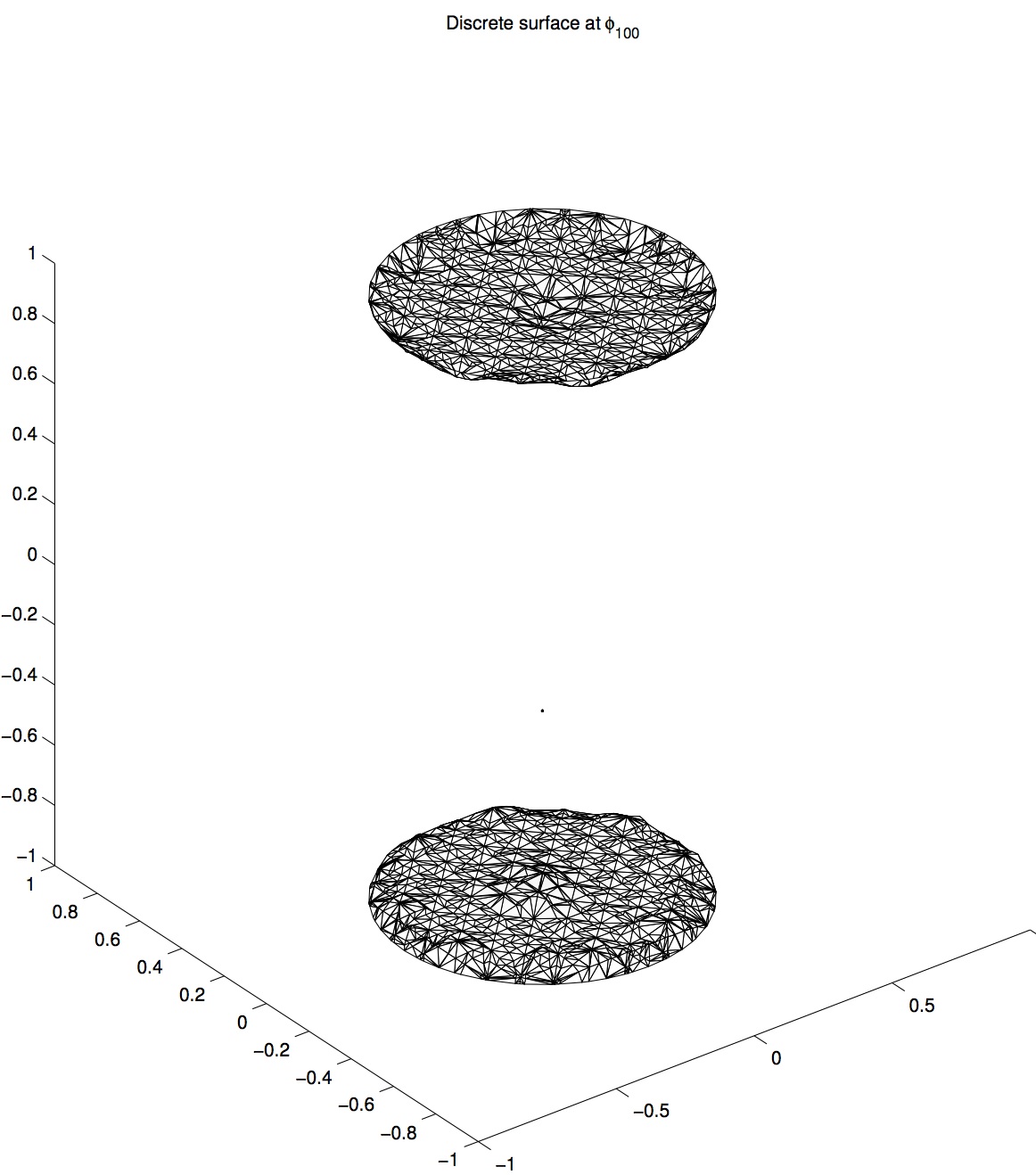}}
\par\end{centering}
\caption{Cut catenoid\label{fig:CollapsingCylinder}}
\end{figure}

\begin{figure}
\begin{centering}
\includegraphics[width=0.5\paperwidth]{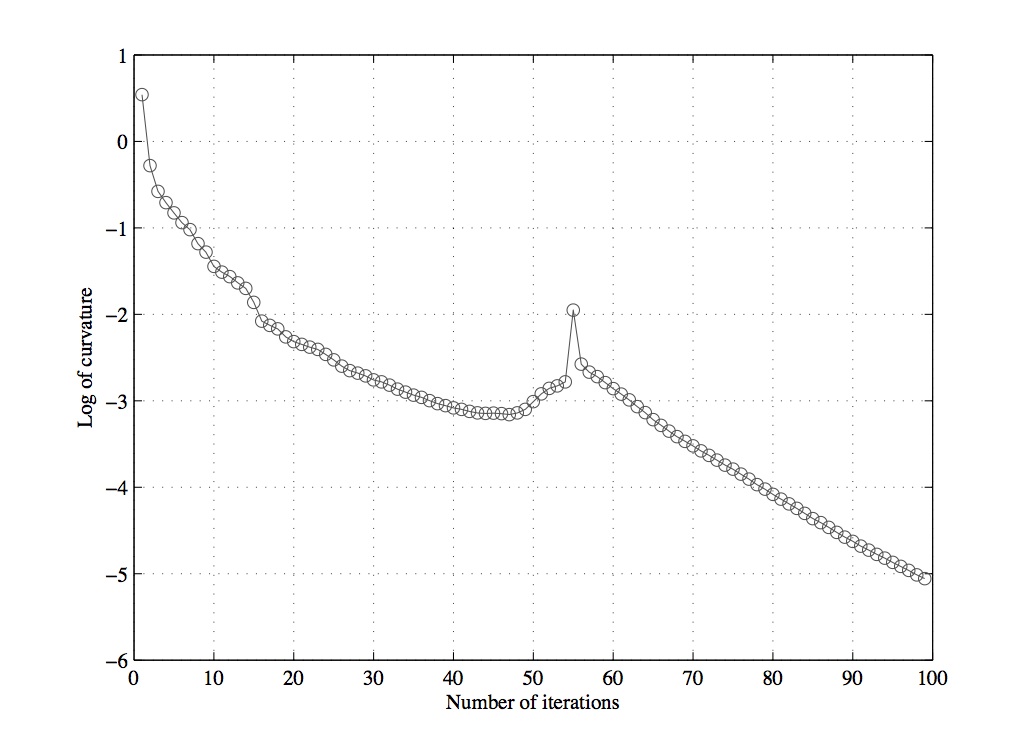}\caption{Convergence of the collapsing cylinder curvature\label{fig:collapsSurfConverg}}
\par\end{centering}
\end{figure}

Fig \ref{fig:CollapsingCylinder} shows a cylinder evolving into two
flat circles. The initial radius is 0.5 and the axis is positioned at $x=0$,
$y=0$. The height is 1 and the mesh domain is of the size $2\times2\times1$.
There are 47017 tetrahedral elements. Figure \ref{fig:collapsSurfConverg}
shows the convergence of mean curvature for the collapsing cylinder.

\begin{figure}
\begin{centering}
\subfloat[Timestep 1]{\includegraphics[width=0.45\columnwidth,height=0.4\paperheight,keepaspectratio]{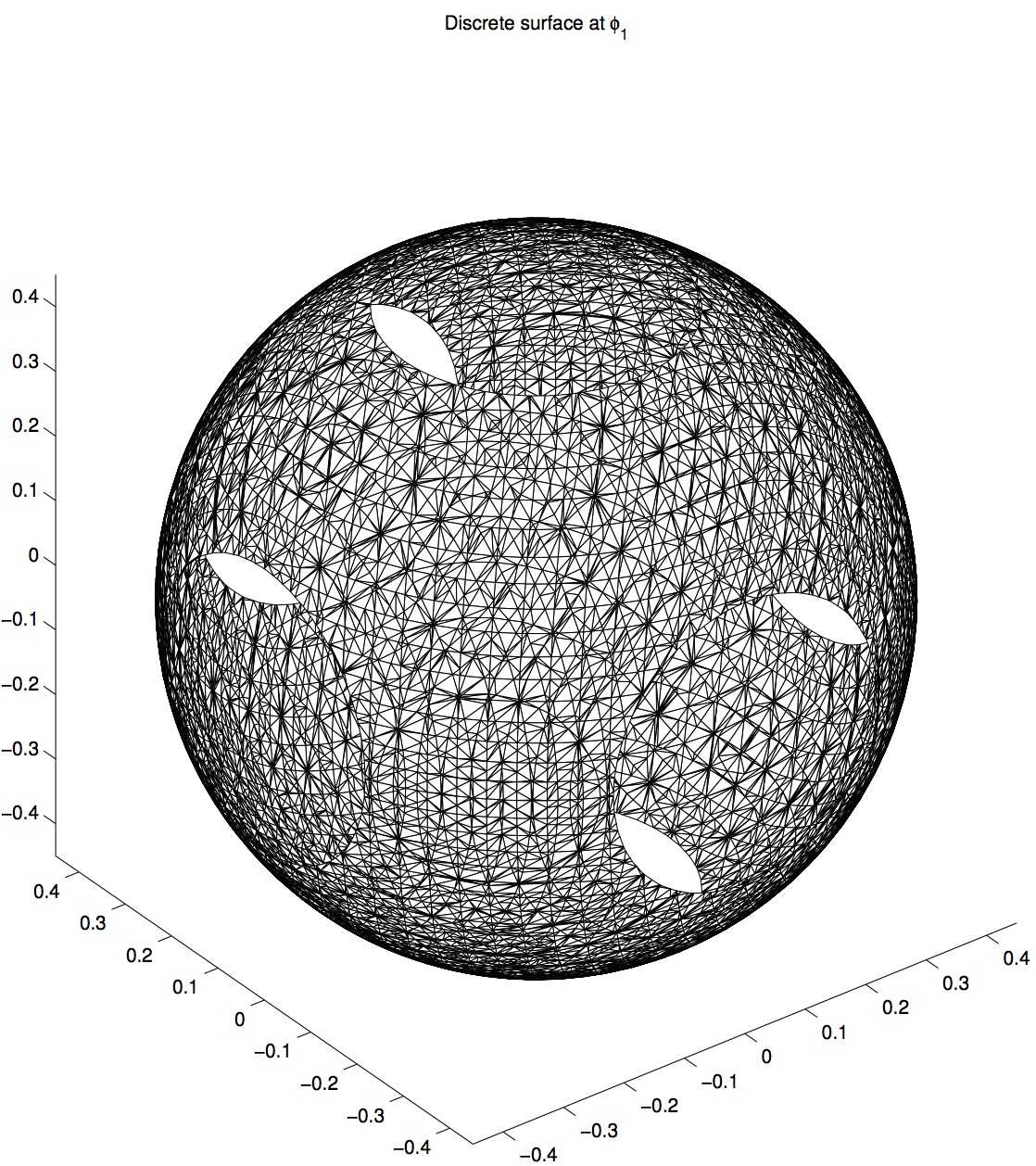}}
\subfloat[Timestep 5]{\includegraphics[width=0.45\columnwidth,height=0.4\paperheight,keepaspectratio]{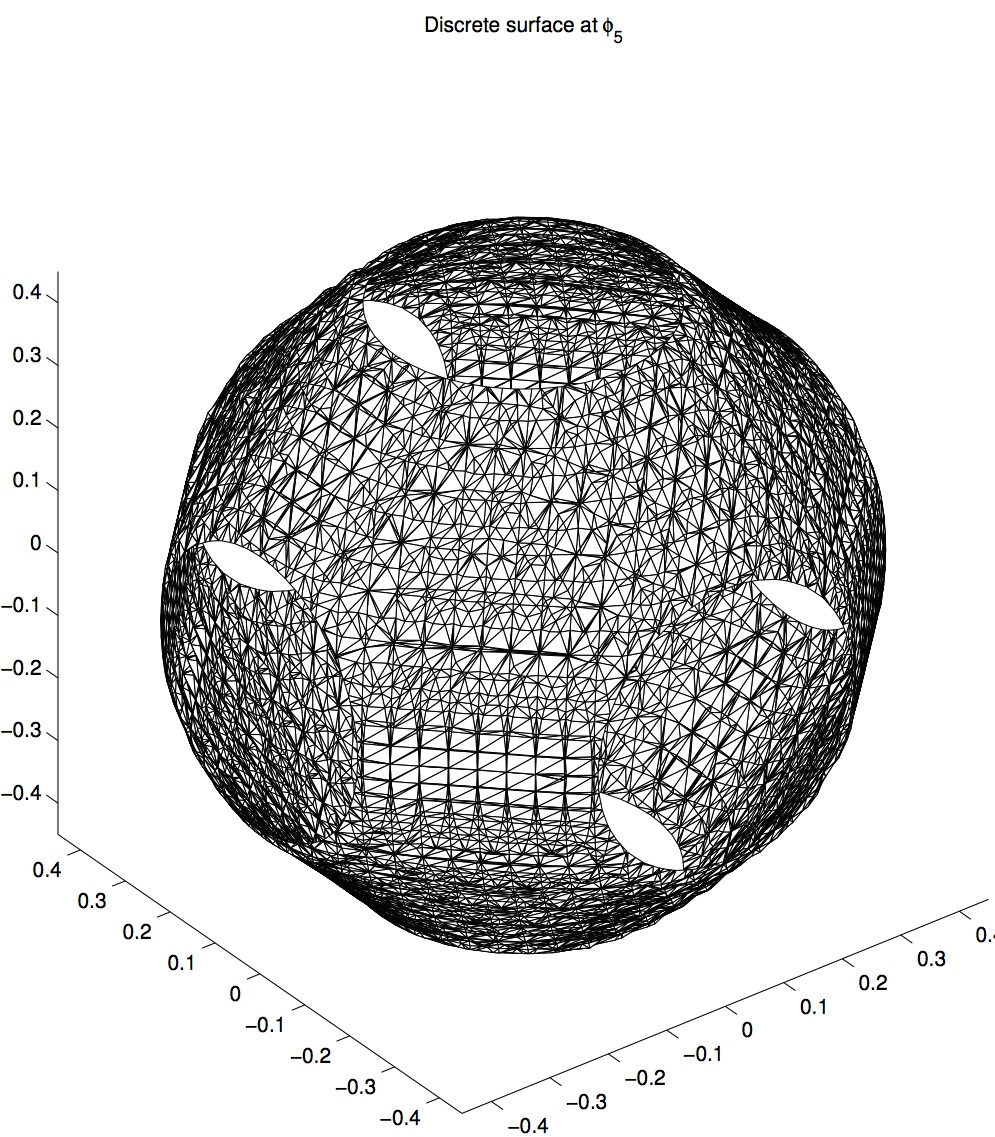}}
\par\end{centering}
\caption{Schwarz minimal surface, early time\label{fig:SchwarzFigure1}}
\end{figure}

\begin{figure}
\begin{centering}
\subfloat[Timestep 10]{\includegraphics[width=0.45\columnwidth,height=0.4\paperheight,keepaspectratio]{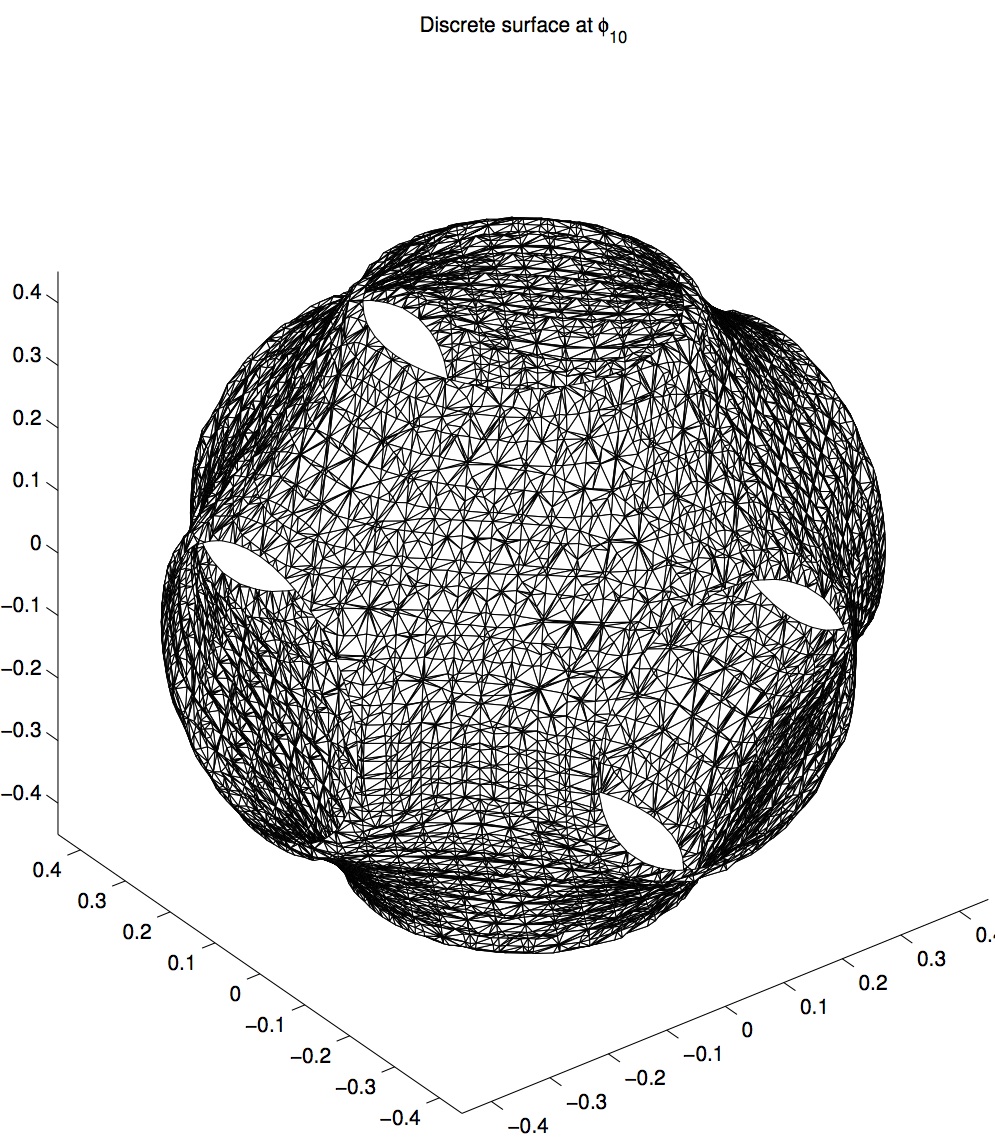}}
\subfloat[Timestep 30]{\includegraphics[width=0.45\columnwidth,height=0.4\paperheight,keepaspectratio]{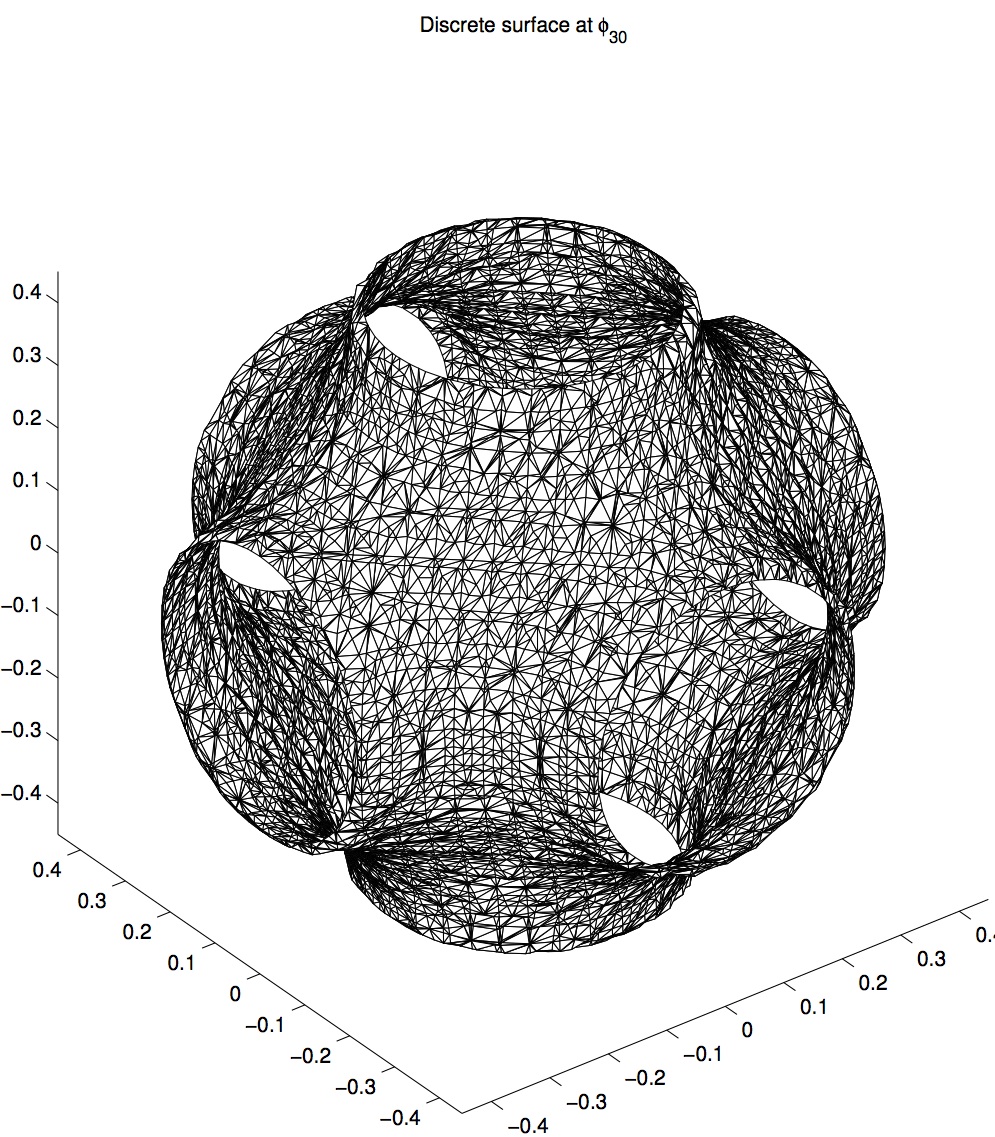}}
\par\end{centering}
\caption{Schwarz minimal surface, intermediate time\label{fig:SchwarzFigure2}}
\end{figure}

\begin{figure}
\begin{centering}
\subfloat[Timestep 60]{\includegraphics[width=0.45\columnwidth,height=0.4\paperheight,keepaspectratio]{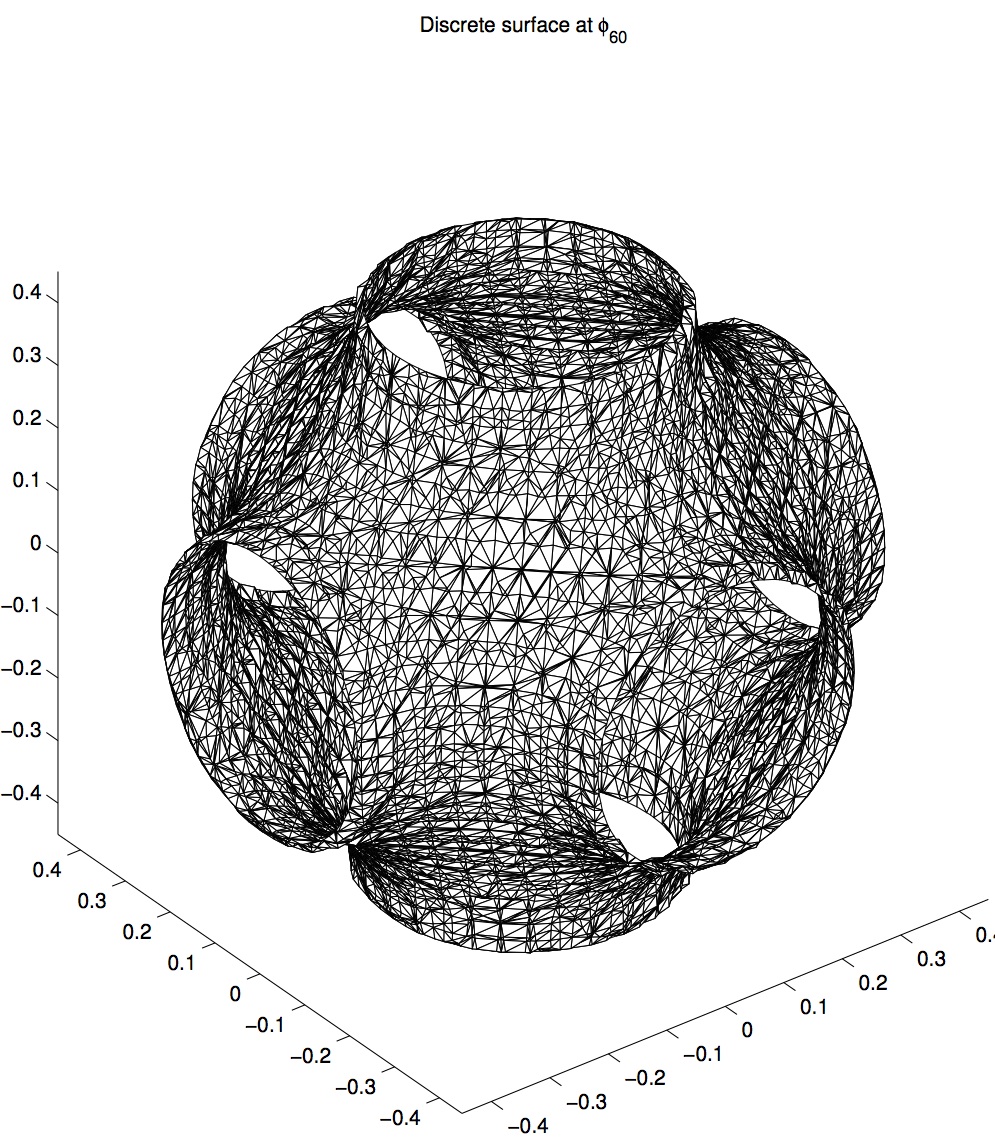}}
\subfloat[Timestep 106]{\includegraphics[width=0.45\columnwidth,height=0.4\paperheight,keepaspectratio]{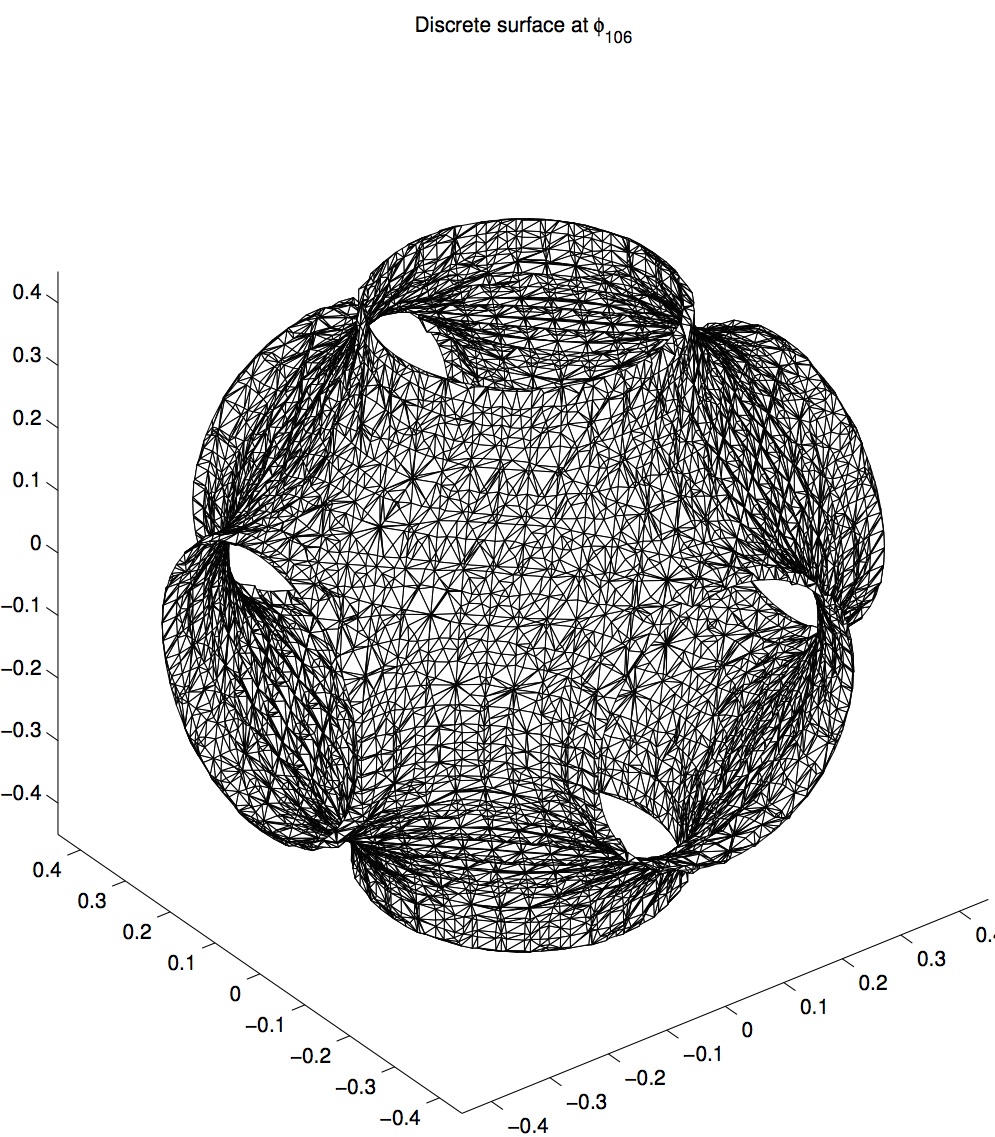}}
\par\end{centering}
\caption{Schwarz minimal surface, late time\label{fig:SchwarzFigure3}}
\end{figure}

\begin{figure}
\includegraphics[width=0.6\paperwidth]{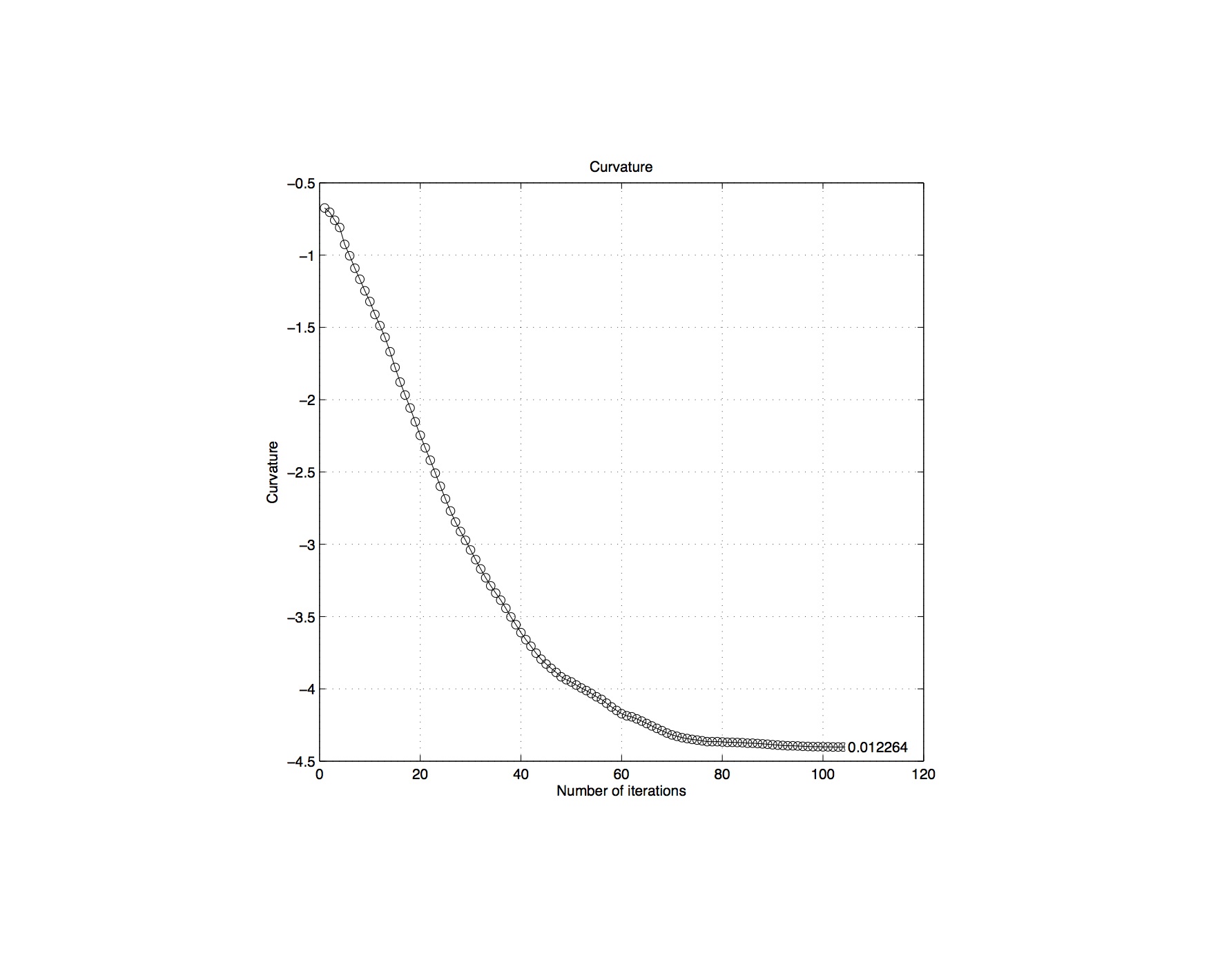}\caption{Convergence of the Schwarz surface\label{fig:SchwarzConverg}}
\end{figure}

Figs \ref{fig:SchwarzFigure1}--\ref{fig:SchwarzFigure3} shows an evolving Schwarz minimal surface starting from a sphere of radius 0.5. The mesh used 196608 tetrahedra of which approximately 24976 were active at any given time. See Figure \ref{fig:SchwarzConverg} for the convergence of the curvature.

\section{Concluding remarks}

We have proposed a novel way of computing minimal surfaces using a 
discretization of the Laplace--Beltrami operator on a 2D surface embedded in a 3D mesh. The approach has a strong theoretical 
foundation in terms of proven accuracy of the computed curvature vector used to drive the evolution of the surface. In future work, 
we will consider more complex surface energies, e.g., leading to membrane elasticity on the surface as in the recent work by Hansbo and Larson  \cite{HaLa14}.


\bibliographystyle{plain}
\addcontentsline{toc}{section}{\refname}\bibliography{minimal}
\newpage

\begin{figure}
\begin{centering}
\includegraphics[width=0.5\paperwidth]{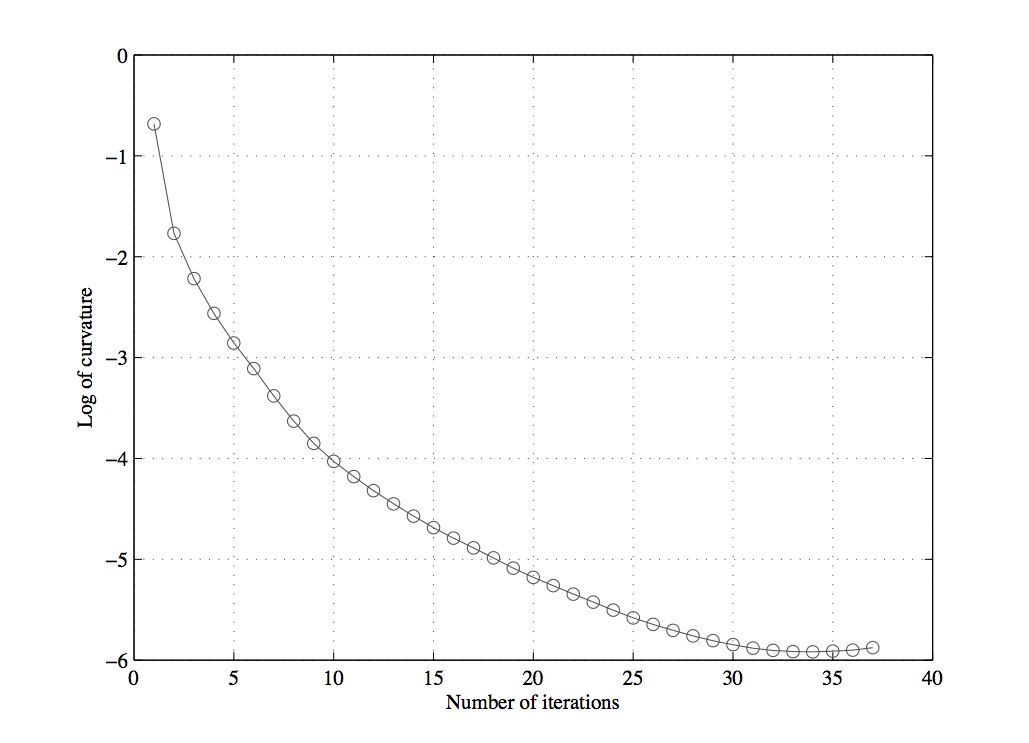}\caption{Convergence of the curvature\label{fig:CatenoidConverg}}

\par\end{centering}

\end{figure}



\end{document}